\documentclass{amsart}

\usepackage{lscape,graphicx}
\usepackage{graphicx}

\usepackage{tikz}
\usetikzlibrary{arrows}
\tikzstyle{block}=[draw opacity=.7,line width=1.4cm]

\theoremstyle{definition}

\theoremstyle{remark}

\numberwithin{equation}{section}

\begin{document}

\title{Fundamental relations between the Dirichlet beta function, Euler numbers, and Riemann zeta function for positive integers}

\author{{Michael A. Idowu}\\
\\
{SIMBIOS Centre, \\
University of Abertay, \\
Dundee DD1 1HG, UK\\
m.idowu@abertay.ac.uk}}

\maketitle

\begin{abstract}
A new definition for the Dirichlet beta function for positive integer arguments is discovered and presented for the first time. This redefinition of the Dirichlet beta function, based on the polygamma function for some special values, provides a general method for obtaining all special constants associated with Dirichlet beta function. We also show various new and fundamental relations between the polygamma function, Riemann zeta, the even-indexed euler numbers, the Dirichlet beta functions in a way never seen or imagined before. 
\end{abstract}


\section{ Introduction}
The polygamma function $\psi^{(k)}(z)$ ought to be given more attention. Following K.S. K\"{o}lbig's suggestion \cite{Kol96, Ido12}, we focus on the special values z=${1\over4}$ and z=${3\over4}$ and obtain another truly remarkable result - a new redefinition for the Dirichlet beta function for positive integer arguments. 

Building on previous work that was based on finding simple methods for deriving the special values of the Riemann zeta function for small integer arguments \cite{Ido12}, we show for the first time how to derive some special values (constants) of the Dirichlet beta function and the even-indexed Euler numbers using either the polygamma $\psi^{(k)}(z)$ or the cotangent function. 

Our results show new relations between the Riemann zeta, Euler numbers, and the Dirichlet beta functions in a way never seen before.  

\subsection{New definition for the Dirichlet beta function}
The well-known Dirichlet beta function is defined as 
\begin{equation}
\beta(s)=\sum_{k=1}^{\infty} {(-1)^{k-1}\over{(2k-1)^s}} = 1-{1 \over {3^s}}+{1 \over {5^s}}-{1 \over {7^s}}+\dots
\end{equation}
; $Re(s)>0$. We wish to demonstrate how the Dirichlet beta function $\beta(s)$ is closely related to the well-known Riemann zeta function
\begin{equation}
\zeta(s)= 1+{1 \over {2^s}}+{1 \over {3^s}}+{1 \over {4^s}}+\dots
\end{equation}
We illustrated in a previous work \cite{Ido12} how the Riemann zeta function may be defined as
\begin{equation}\label{zeta}
\zeta(s)= (-1)^s.{{( \psi^{(s-1)}({1\over4})+\psi^{(s-1)}({3\over4}))}\over{2^s.(2^s-1)}}.{1 \over{\Gamma(s)}}
\end{equation}
for any positive integer argument s, where the polygamma function $\psi^{(s-1)}(x)$ is defined as 
\begin{equation}
\psi^{(s-1)}(x) = {d^{s-1}\over{dx^{s-1}}}\psi(x) = {d^{s}\over{dx^{s}}}\ln \Gamma(x).
\end{equation}
Further manipulation of the polygamma function reveals that the Dirichlet beta function may be represented as 
\begin{equation}\label{beta}
\beta(s)= (-1)^{s}({{( \psi^{(s-1)}({1\over4})-\psi^{(s-1)}({3\over4}))}\over{2^s.2^s}}){1 \over{\Gamma(s)}}=\sum_{k=1}^{\infty} {(-1)^{k-1}\over{(2k-1)^s}}.
\end{equation}

\subsection{Derivation of a new definition for the Dirichlet beta function}
From the results in $\ref{zeta}$ we infer that
\begin{equation}\label{Oddzeta}
{{2^s-1}\over{2^s}}\zeta(s)= (-1)^s.{{( \psi^{(s-1)}({1\over4})+\psi^{(s-1)}({3\over4}))}\over{2^s.2^s}}.{1 \over{\Gamma(s)}}
=\sum_{k=1}^{\infty} {1\over{(2k-1)^s}} 
\end{equation}
Relating $\ref{beta}$ to $\ref{Oddzeta}$ the following relation emerges:
\begin{equation}\label{beta2}
\beta(s)= {{2^s-1}\over{2^s}}\zeta(s)-2\sum_{k=1}^{\infty} {1\over{(4k-1)^s}} ={{2^s-1}\over{2^s}}\zeta(s)-{{2}\over{2^s.2^s}}\sum_{k=0}^{\infty} {1\over{(k+{3\over4})^s}} 
\end{equation}
Therefore,
\begin{equation}\label{beta3}
\beta(s)= {{2^s-1}\over{2^s}}\zeta(s)-(-1)^s{{2}\over{2^s.2^s}}{1 \over{\Gamma(s)}}{{ \psi^{(s-1)}({3\over4})}}.
\end{equation}
Hence, the result in \ref{beta} is inferred from the polygamma relation
\begin{equation}\label{beta4}
\beta(s)= {{2^s-1}\over{2^s}}(-1)^s{{( \psi^{(s-1)}({1\over4})+\psi^{(s-1)}({3\over4}))}\over{2^s.(2^s-1)}}{1 \over{\Gamma(s)}}-(-1)^s{{2}\over{2^s.2^s}}{1 \over{\Gamma(s)}}{{ \psi^{(s-1)}({3\over4})}}
\end{equation}
that is obtained after instant substitution. 

\subsection{Derivation of new definitions for the Riemann zeta function}
From equation \ref{beta3}, variant definitions for the Riemann zeta function may be derived (in terms of the Dirichlet beta function), i.e. 
\begin{equation}\label{beta3}
\zeta(s)={{2^s}\over{2^s-1}}\beta(s)+(-1)^s{{2}\over{2^s(2^s-1)}}{1 \over{\Gamma(s)}}{{ \psi^{(s-1)}({3\over4})}}.
\end{equation}
This implies that for even integer arguments, 
\begin{equation}\label{beta3}
\zeta(2s)={{2^{2s}}\over{2^{2s}-1}}\beta(2s)+{{2}\over{2^{2s}(2^{2s}-1)}}{1 \over{\Gamma(2s)}}{{ \psi^{(2s-1)}({3\over4})}},
\end{equation}
and for odd integer arguments \footnote{It is the author's original view that such dual-nature implicit in Riemann zeta (re)definitions ought to be given more attention if we are ever going to find exact expressions for $\zeta(2s+1)$ in terms of other well-known mathematical constants. \copyright 2012 Michael Idowu.}, 
\begin{equation}\label{beta3}
\zeta(2s+1)={{2^{2s+1}}\over{2^{2s+1}-1}}\beta(2s+1)-{{2}\over{2^{2s+1}(2^{2s+1}-1)}}{1 \over{\Gamma(2s+1)}}{{ \psi^{(2s)}({3\over4})}}
\end{equation}, where $s>0$ is a positive integer number.

\section{Elegant Expression for the Dirichlet beta function} 
The reflection relation 
\begin{equation}
(-1)^{s-1} \psi^{(s-1)}({1-z})-\psi^{(s-1)}({z}) = \pi{d^{(s-1)}\over{dz^{(s-1)}}}cot(\pi z) 
\end{equation}
may be used to provide alternative definition for the Dirichlet beta function, i.e.
\begin{equation}\label{zetaCot}
- \psi^{(2s-1)}({1-z})-\psi^{(2s-1)}({z}) = \pi{d^{(2s-1)}\over{dz^{(2s-1)}}}cot(\pi z) 
\end{equation}
\begin{equation}\label{betaCot}
 \psi^{(2s)}({1-z})-\psi^{(2s)}({z}) = \pi{d^{(2s)}\over{dz^{(2s)}}}cot(\pi z) 
\end{equation}
As mentioned in a previous work \cite{Ido12}, an equation such as \ref{zetaCot} is useful for composing the actual values of the Riemann zeta function for even integer arguments. Since equation \ref{betaCot} is more consistent with our new definition of $\beta(s)$ provided in \ref{beta}, we employ the relation in \ref{betaCot} to compose the actual values of $\beta{2s+1}$. 

In other words, expressing (or composing) the exact values of $\beta(2s+1)$ is a simple as composing the values of $\zeta(2s)$ as demonstrated in our previous work \cite{Ido12}.  Also it might be worth mentioning that finding exact expressions for $\beta(2s)$ is as hard as finding expressions for $\zeta(2s+1)$ as indicated in by K.S. K\"{o}lbig in his paper \cite{Kol96}many years ago. 

To demonstrate the usefulness of the relation specified in \ref{betaCot}, we employ the equivalent relation 
\begin{equation}
{2^{2s+1}(2^{2s+1})\Gamma(2s+1)}\beta(2s+1) = \pi {d^{(2s)}\over{dz^{(2s)}}}cot(\pi z) \mid{_{z\rightarrow{1\over4}}}
\end{equation} 
and derive some well-known constants associated with $\beta(2s+1)$, e.g. the actual values of $\beta(1),\beta(3), \beta(5), \beta(7).$

\subsection{Representation of some special constants of the Dirichlet beta function}
Some well-known constants of $\beta(s)$ may be calculated as follows:
\[
\beta(1)= (-1)^1.{{( \psi({1\over4})-\psi({3\over4}))}\over{2^1.(2^1).\Gamma(1)}} = {\pi \over 4};
\]
\[
\beta(2)= (-1)^2.{{( \psi^{'}({1\over4})-\psi^{'}({3\over4}))}\over{2^2.(2^2).\Gamma(2)}} = G ,
\]
where G represents the Catalan's constant $\approx 0.915965594177$;
\[
\beta(3)= (-1)^3.{{( \psi^{''}({1\over4})-\psi^{''}({3\over4}))}\over{2^3.(2^3).\Gamma(3)}} \approx 0.968946146259;
\]
\[
\beta(4)= (-1)^4.{{( \psi^{'''}({1\over4})-\psi^{'''}({3\over4}))}\over{2^4.(2^4).\Gamma(4)}} \approx 0.98894455174; 
\]
\[
\beta(5)= (-1)^5.{{( \psi^{''''}({1\over4})-\psi^{''''}({3\over4}))}\over{2^5.(2^5).\Gamma(5)}} \approx 0.996157828077; 
\]

\subsection{General formula for the Dirichlet beta function for odd integers}\label{betaSect}
Using the formula 
\begin{equation}\label{betaF}
\beta(2s+1) ={{ \pi {d^{(2s)}\over{dz^{(2s)}}}cot(\pi z) \mid{_{z\rightarrow{1\over4}}}} \over 
                        {{2^{2s+1}(2^{2s+1})\Gamma(2s+1)}} },
\end{equation}
we derive
\small
\[
\beta(1) ={{ \pi cot(\pi z) \mid{_{z\rightarrow{1\over4}}}} \over 
                        {{2^{1}(2^{1})\Gamma(1)}} } = {{\bf{1}}\pi \over {{2(2^{1})\Gamma(1)}}},
\]
\[
\beta(3) ={{ \pi {d^{(2)}\over{dz^{(2)}}}cot(\pi z) \mid{_{z\rightarrow{1\over4}}}} \over 
                        {{2^{3}(2^{3})\Gamma(3)}} } = 
{{ \pi (2\pi^2)cot(\pi z)(cot(\pi z)^2+1) \mid{_{z\rightarrow{1\over4}}}} \over 
                        {{2^{3}(2^{3})\Gamma(3)}} } = { {2\pi^3(1)(1+1)}\over 
                        {{2^{3}(2^{3})\Gamma(3)}} } = {{\bf{1}}\pi^3 \over {{2(2^{3})\Gamma(3)}}},
\]
\[
\beta(5) ={{ \pi {d^{(4)}\over{dz^{(4)}}}cot(\pi z) \mid{_{z\rightarrow{1\over4}}}} \over 
                        {{2^{5}(2^{5})\Gamma(5)}} } = 
{{ \pi ((16\pi^4)(1)(1+1)^2+(8\pi^4)(1)^3(1+1)) \mid{_{z\rightarrow{1\over4}}}} \over 
                        {{2^{5}(2^{5})\Gamma(5)}} } = { {(64+16)\pi^5}\over 
                        {{2^{5}(2^{5})\Gamma(5)}} } = {{\bf{5}}\pi^5 \over {{2(2^{5})\Gamma(5)}}},
\]
\[
\beta(7) ={{ \pi {d^{(6)}\over{dz^{(6)}}}cot(\pi z) \mid{_{z\rightarrow{1\over4}}}} \over 
                        {{2^{7}(2^{7})\Gamma(7)}} } = 
{{ \pi (272\pi^6.1.(1+1)^3+32\pi^6.1^5(1+1)+416\pi^6.1^3(1+1)^2) \mid{_{z\rightarrow{1\over4}}}} \over 
                        {{2^{7}(2^{7})\Gamma(7)}} } = { {\bf{61}\pi^7}\over 
                        {{2(2^{7})\Gamma(7)}} }.
\]

\normalsize

\subsection{General formula for the even-indexed Euler numbers}
Due to the emergence of the Euler numbers observed in section \ref{betaSect}, we boldly guess and confirm that the next value to be derived is: 
\[
\beta(9) = { {\mid {\bf{E_8}} \mid\pi^9}\over {{2(2^{9})\Gamma(9)}} } ={ {\bf{1385}\pi^9}\over {{2(2^{9})\Gamma(9)}} },
\]where ${E_8}$ is the related Euler number for the observed pattern.
Isn't this last derivation elegant? Even more elegant is the following bolder guess that the general formula for all even-indexed Euler numbers is
\begin{equation}
E_{2s} =  {{(2^{2s+1})\Gamma(2s+1)\beta(2s+1)}\over{(\pi i)^{2s+1}}}.2i
\end{equation}

\section{Summary of Main Results}
To summarise the results obtained so far
\begin{equation}\label{zetaF}
\zeta(2s) = (-1)^{2s}({{( \psi^{(2s-1)}({1\over4})+\psi^{(2s-1)}({3\over4}))}\over{2^{2s}(2^{2s}-1)}}){1 \over{\Gamma(2s)}} = {{ \pi {d^{(2s-1)}\over{dz^{(2s-1)}}}cot(\pi z) \mid{_{z\rightarrow{1\over4}}}} \over 
                        {{2^{2s}(2^{2s}-1)\Gamma(2s)}} }
\end{equation}

\begin{equation}\label{betaF}
\beta(2s+1) =(-1)^{2s+1}({{( \psi^{(2s)}({1\over4})-\psi^{(2s)}({3\over4}))}\over{2^{2s+1}(2^{2s+1})}}){1 \over{\Gamma(2s+1)}}={{ \pi {d^{(2s)}\over{dz^{(2s)}}}cot(\pi z) \mid{_{z\rightarrow{1\over4}}}} \over 
                        {{2^{2s+1}(2^{2s+1})\Gamma(2s+1)}} },
\end{equation}

\begin{equation}\label{EulerF}
 E_{2s} =(-1)^{2s+1}({{( \psi^{(2s)}({1\over4})-\psi^{(2s)}({3\over4}))}\over{(2\pi i)^{2s+1}}}).2i = {{ \pi {d^{(2s)}\over{dz^{(2s)}}}cot(\pi z) \mid{_{z\rightarrow{1\over4}}}} \over {(2\pi i)^{2s+1}}}.2i
\end{equation}

These new reformulations enable a new conceptualisation of the various close relations between the Euler numbers, Dirichlet beta function, and Riemann zeta function for positive integer numbers.

\end{document}